\documentclass[10pt,a4paper]{amsart}

\usepackage{amssymb}
\usepackage{amsthm}
\usepackage{url}
\usepackage{amsfonts}
\usepackage{mathrsfs}
\usepackage{graphicx}
\usepackage{color}
\usepackage{latexsym}
\usepackage{rotating}
\usepackage{xspace}
\usepackage[all]{xy}
\usepackage{longtable}
\usepackage{leftidx}
\usepackage{mathtools}
\usepackage{titlesec}
\usepackage{tikz}
\usetikzlibrary{calc}
\usetikzlibrary{arrows}
\usetikzlibrary{decorations.markings}
\usepackage{geometry}

\newtheorem{theorem}{Theorem}[section]
\numberwithin{equation}{theorem}
\newtheorem{lemma}[theorem]{Lemma}
\newtheorem{proposition}[theorem]{Proposition}
\newtheorem{corollary}[theorem]{Corollary}

\theoremstyle{definition}
\newtheorem{definition}[theorem]{Definition}
\theoremstyle{fact}
\newtheorem{fact}[theorem]{Fact}
\newtheorem{example}[theorem]{Example}

\newtheorem{remark}[theorem]{Remark}

\theoremstyle{conjecture}

\newcommand{\Ass}{\operatorname{Ass}}
\newcommand{\f}{\operatorname{f}}

\newcommand{\grade}{\operatorname{grade}}

\newcommand{\ara}{\operatorname{ara}}

\newcommand{\cd}{\operatorname{cd}}

\newcommand{\mcd}{\operatorname{mcd}}

\newcommand{\Ext}{\operatorname{Ext}}

\newcommand{\Ann}{\operatorname{Ann}}

\newcommand{\Rad}{\operatorname{Rad}}

\newcommand{\depth}{\operatorname{depth}}

\newcommand{\fm}{\frak{m}}
\newcommand{\fp}{\frak{p}}
\newcommand{\fq}{\frak{q}}
\newcommand{\fa}{\frak{a}}

\newcommand{\fn}{\frak{n}}

\newcommand{\suchthat}{\;\ifnum\currentgrouptype=16 \middle\fi|\;}

\newenvironment{prf}[1][Proof]{\begin{proof}[\bf #1]}{\end{proof}}

\makeatletter
\newcommand{\holim@}[2]{%
\vtop{\m@th\ialign{##\cr
\hfil$#1\operator@font holim$\hfil\cr
\noalign{\nointerlineskip\kern1.5\ex@}#2\cr
\noalign{\nointerlineskip\kern-\ex@}\cr}}%
}
\newcommand{\holim}{%
\mathop{\mathpalette\holim@{\rightarrowfill@\textstyle}}\nmlimits@
}
\makeatother

\makeatletter
\def\@secnumfont{\bfseries}
\def\section{\@startsection{section}{1}%
\z@{.7\linespacing\@plus\linespacing}{.5\linespacing}%
{\normalfont\Large\bfseries\filcenter}}
\def\subsection{\@startsection{subsection}{2}%
\z@{.5\linespacing\@plus.7\linespacing}{-.5em}%
{\normalfont\large\bfseries}}
\makeatother
\begin{document}

\author[Ahad Rahimi]
{Ahad Rahimi}

\title[Homological Properties of $L\otimes_{\Bbbk} N$]
{Homological Properties of $L\otimes_{\Bbbk} N$}

\address{A. Rahimi, Department of Mathematics, Razi University, Kermanshah, Iran.}
\email{ahad.rahimi@razi.ac.ir}



\begin{abstract}
Let $A$ and $B$ be commutative Noetherian algebras over an arbitrary field $\Bbbk$ such that $A\otimes_\Bbbk B$ is Noetherian.
 Let $I$ and $J$ be ideals of $A$ and $B$, respectively, and let $L$ and $N$ be nonzero finitely generated modules over $A$ and $B$, respectively.
  In \cite{R4}, the author investigates the relationships between certain algebraic properties of $L\otimes_\Bbbk N$ and
   the corresponding properties of $L$ and $N$, and vice versa. In the present paper, we characterize the Gorensteinness,
    approximate Cohen--Macaulayness, and maximal depth property of $L\otimes_\Bbbk N$ with respect
     to $\mathcal{I}=I\otimes_\Bbbk B+A\otimes_\Bbbk J$ in terms of the corresponding properties of $L$ and $N$ with respect to $I$ and $J$, respectively.
     In particular, when $A$ and $B$ are local rings, we describe these properties of $L \otimes_{\Bbbk} N$ as a module over the local ring $A \otimes_\Bbbk B$.
\end{abstract}

\maketitle

\tableofcontents

\section{Introduction}
Throughout this article, $R$ denotes a commutative Noetherian ring with identity.

Let $A$ and $B$ be commutative Noetherian $\Bbbk$-algebras over a field
$\Bbbk$, and assume that $A\otimes_\Bbbk B$ is also Noetherian.
Let $I$ and $J$ be ideals of $A$ and $B$, respectively, and set
$
\mathcal{I}=I\otimes_\Bbbk B+A\otimes_\Bbbk J.
$

We may also consider the case in which $A\otimes_\Bbbk B$ is a local
ring with maximal ideal
$
\mathfrak{M}
=\mathfrak{m}\otimes_\Bbbk B+A\otimes_\Bbbk\mathfrak{n},
$
where $(A,\mathfrak{m})$ and $(B,\mathfrak{n})$ are local
$\Bbbk$-algebras. In this setting, $A$ and $B$ have $\Bbbk$ as their
common residue field, and at least one of $A$ and $B$ is algebraic
over $\Bbbk$.
Furthermore, let $L$ and $N$ be nonzero finitely generated modules
over $A$ and $B$, respectively(see \cite{R4}).

 In this paper, we characterize the
Gorenstein property, approximate Cohen--Macaulayness, and the maximal
depth property of $L\otimes_\Bbbk N$ with respect to $\mathcal{I}$
in terms of the corresponding properties of $L$ and $N$ with respect
to $I$ and $J$, respectively. We also answer affirmatively a question raised in \cite{R5}:
 if $M=L\otimes_\Bbbk N$ satisfies certain homological properties with respect to $A\otimes_\Bbbk J$,
  then $M/(I\otimes_\Bbbk B)M$ inherits the corresponding properties with respect to $\mathcal{I}$.
   However, Gorensteinness and maximal Cohen--Macaulayness are not preserved in general under this quotient construction.

We begin in Section 2 by investigating the Gorensteinness of
$L\otimes_{\Bbbk} N$.  We say that an $R$-module $K$
is Gorenstein with respect to an ideal $\mathfrak{a}$ of $R$ if
$\Ext^i_R(R/\mathfrak{a},K)$ is nonzero for exactly one value of
$i$, namely $i=\cd(\mathfrak{a},R)$ (see \cite{PDR}).
In the usual local setting $(R,\mathfrak{m'})$, a Gorenstein
$R$-module is one that is Gorenstein with respect to $\mathfrak{m'}$.
Using a result from \cite{R4} together with the Künneth formula,
we prove that $L\otimes_{\Bbbk} N$ is Gorenstein with respect to
$\mathcal{I}$ if and only if $L$ and $N$ are Gorenstein with respect
to $I$ and $J$, respectively; see Theorem~\ref{2.4}.
In particular, in the local setting, $L\otimes_{\Bbbk} N$ is
Gorenstein if and only if both $L$ and $N$ are Gorenstein.

We next establish in Lemma \ref{3.2} that if $\Bbbk$ is an algebraically closed field, then
  $\Rad(\mathcal{I})= \Rad(I) \otimes_\Bbbk B + A \otimes_\Bbbk \Rad(J)$.  Using this fact, we prove in Theorem~\ref{3.4} that
 \[
\ara(\mathcal{I}, L \otimes_\Bbbk N)
= \ara(I, L) + \ara(J, N).
\]
This result is used to establish the main result of the next section.

We next examine the approximate Cohen--Macaulayness of $L\otimes_\Bbbk N$.
By a result of \cite{FR}, approximately
Cohen--Macaulay modules are precisely those that are sequentially
Cohen--Macaulay and whose Cohen--Macaulay filtrations have length
two. Thus, approximate Cohen--Macaulayness can be viewed as a
natural extension of Cohen--Macaulayness.
Example \ref{4.5} shows that the approximate Cohen--Macaulayness of $L$ and $N$ does
 not imply the approximate Cohen--Macaulayness of $L\otimes_\Bbbk N$. Using Theorem~\ref{3.4}, together with the assumption
that $\Bbbk$ is algebraically closed and the equalities
$\cd(I,L)=\ara(I,L)$ and $\cd(J,N)=\ara(J,N)$, we prove the
following characterization: $L\otimes_\Bbbk N$ is approximately
Cohen--Macaulay with respect to $\mathcal{I}$ but not
Cohen--Macaulay with respect to $\mathcal{I}$ if and only if
$L$ is approximately Cohen--Macaulay with respect to $I$ but not
Cohen--Macaulay with respect to $I$, while $N$ is
Cohen--Macaulay with respect to $J$; see Theorem~\ref{4.7}.

In particular, in the local setting, $L\otimes_\Bbbk N$ is
approximately Cohen--Macaulay but not Cohen--Macaulay if and only
if $L$ is approximately Cohen--Macaulay but not Cohen--Macaulay,
and $N$ is Cohen--Macaulay.

Our next result concerns the maximal depth property, which
generalizes approximate Cohen--Macaulayness. An $R$-module $K$ is
said to have maximal depth with respect to an ideal $\mathfrak{a}$
if
$
\grade(\mathfrak{a},K)=\mcd(\mathfrak{a},K),
$
where
$
\mcd(\mathfrak{a},K)
=\min\{\cd(\mathfrak{a},R/\mathfrak{p}):
\mathfrak{p}\in\Ass_R K\}.
$
We first show that, if $\Bbbk$ is an algebraically closed field, then
\[
\mcd(\mathcal{I},L\otimes_\Bbbk N)
=\mcd(I,L)+\mcd(J,N);
\]
see Proposition~\ref{5.1}. Using this equality, we prove that
$L\otimes_\Bbbk N$ has maximal depth with respect to $\mathcal{I}$
if and only if $L$ and $N$ have maximal depth with respect to $I$
and $J$, respectively; see Theorem~\ref{5.2}. In particular, in
the local setting, $L\otimes_\Bbbk N$ has maximal depth if and only
if both $L$ and $N$ have maximal depth.

In the final section, we investigate similar questions raised in \cite{R5} (see also \cite{NR})
 concerning $L\otimes_\Bbbk N$. More precisely, we show that if $M=L\otimes_\Bbbk N$ is Cohen--Macaulay,
  generalized Cohen--Macaulay, approximately Cohen--Macaulay, or sequentially Cohen--Macaulay with
   respect to $A\otimes_\Bbbk J$, then $M/(I\otimes_\Bbbk B)M$ has the corresponding property with
    respect to $\mathcal{I}$. Moreover, if $M$ has maximal depth with respect to $A\otimes_\Bbbk J$,
     then so does $M/(I\otimes_\Bbbk B)M$ with respect to $\mathcal{I}$. In contrast, Gorensteinness
      and maximal Cohen--Macaulayness are not preserved in general under this quotient construction; see Remarks~\ref{6.10} and~\ref{6.11}.
As an application, we recover Propositions 6.1--6.3 of \cite{R5}.

\section{Gorensteiness of $L\otimes_{\Bbbk} N$}
Let $\fa$ be an ideal of $R$, and let $K$ be a finitely generated $R$-module.
 The cohomological dimension of $K$ with respect to $\fa$, denoted by $\cd(\fa, K)$, is defined as the supremum of
all integers $i$ for which $H^i_\fa(K)\neq 0$, where $H^i_\fa(K)$ is the $i$-th local cohomology module of $K$ with respect to $\fa$.
 It is important to note that, $K=\fa K$ if and only if $\cd(\fa, K)=-\infty$.  The grade of $K$ with respect to $\fa$, denoted by $\grade(\fa, K)$,
 is defined as the infimum of all integers $i$ such that $H^i_\fa(K)\neq 0$.
 It is worth mentioning that $K=\fa K$ if and only if $\grade(\fa, K)=+\infty$. Note that if $\fa K\neq K$, then $\grade(\fa, K)\leq \cd(\fa, K)$.
\begin{definition}
\label{2.1}
We say that $K$ is \emph{Gorenstein with respect to $\fa$} if $\Ext_R^i(R/\fa, K)=0$ for all $i\neq \cd(\fa,R)$ (see \cite{PDR}).
\end{definition}
Let $(R,\fm')$ be a local ring. Then $K$ is Gorenstein if and only if it is Gorenstein with respect to $\fm'$.

Let $\Bbbk$ be a field, and let $A$ and $B$ be Noetherian $\Bbbk$-algebras such that $A\otimes_{\Bbbk} B$ is Noetherian. Let $I$ and $J$ be ideals of $A$ and $B$, respectively, and set
$
\mathcal{I}=I\otimes_{\Bbbk} B+A\otimes_{\Bbbk} J.
$
Furthermore, let $L$ and $N$ be nonzero finitely generated modules over $A$ and $B$, respectively.

In this section, we characterize the Gorensteinness of $L\otimes_{\Bbbk} N$ with respect to $\mathcal{I}$ in terms of the Gorensteinness of $L$ and $N$ with respect to $I$ and $J$, respectively.
We begin by recalling some facts that will be used in the sequel.

 \begin{fact}{\em
\label{2.2}
In view of \cite[Theorem 2.6]{R4}, we have
\begin{itemize}
\item[{(a)}] $\grade(\mathcal{I}, L\otimes_{\Bbbk} N)=\grade(I, L)+\grade(J, N)$.
\item[{(b)}] $\cd(\mathcal{I}, L\otimes_{\Bbbk} N)=\cd(I, L)+\cd(J, N)$.
\end{itemize}
 }
\end{fact}

We have the following fact
 \begin{eqnarray}	
\label{2.2.1}
(A\otimes_{\Bbbk} B)/\mathcal{I}\cong A/I\otimes_\Bbbk B/J.
\end{eqnarray}

We recall a consequence of the K\"unneth formula.
\begin{fact}
\label{2.3}
\em
Let $\Bbbk$ be a field, and let $A$ and $B$ be $\Bbbk$-algebras. Let $L$ and $L'$ be $A$-modules, with $L$ finitely generated, and let $N$ and $N'$ be $B$-modules, with $N$ finitely generated. Then, for each $i\geq 0$, there is an isomorphism of $A\otimes_{\Bbbk}B$-modules
\[
\Ext^i_{A\otimes_{\Bbbk}B}
(L\otimes_{\Bbbk}N,L'\otimes_{\Bbbk}N')
\cong
\bigoplus_{s+t=i}
\Ext^s_A(L,L')\otimes_{\Bbbk}
\Ext^t_B(N,N').
\]
See, for example, \cite[Proposition A.1.5]{W}.
\end{fact}

\begin{theorem}
\label{2.4}
Continue with the notation and assumptions above. Then the following
conditions are equivalent:
\begin{itemize}
\item[(a)] $L\otimes_{\Bbbk} N$ is Gorenstein with respect to
$\mathcal{I}$;
\item[(b)] $L$ and $N$ are Gorenstein with respect to $I$ and $J$,
respectively.
\end{itemize}
\end{theorem}

\begin{proof}
$(\mathrm{a})\Rightarrow(\mathrm{b})$:
Suppose that $L\otimes_{\Bbbk} N$ is Gorenstein with respect to
$\mathcal{I}$. Then
\[
\Ext^i_{A\otimes_{\Bbbk}B}
\left((A\otimes_{\Bbbk}B)/\mathcal{I},L\otimes_{\Bbbk}N\right)=0
\quad\text{for all}\quad
i\neq \cd(\mathcal{I},A\otimes_{\Bbbk}B).
\]
 Consequently,
\[
\grade(\mathcal{I},L\otimes_{\Bbbk}N)
=\cd(\mathcal{I},A\otimes_{\Bbbk}B)
=\cd(I,A)+\cd(J,B).
\]
 The second equality follows from Fact~\ref{2.2}(b). We also have
$
\grade(\mathcal{I},L\otimes_{\Bbbk}N)
=\grade(I,L)+\grade(J,N)
$
by Fact~\ref{2.2}(a).
Hence
\[
\grade(I,L)+\grade(J,N)
=\cd(I,A)+\cd(J,B).
\]
Since  $\grade(I,L)\leq\cd(I,A)$ and $\grade(J,N)\leq\cd(J,B)$,
we obtain
\[
\grade(I,L)=\cd(I,A)
\quad\text{and}\quad
\grade(J,N)=\cd(J,B).
\]
It follows that
\[
\Ext^s_A(A/I,L)=0
\quad\text{for all }s<\cd(I,A) \  \  \   (*),
\]
and
\[
\Ext^t_B(B/J,N)=0
\quad\text{for all }t<\cd(J,B) \  \  \   (**).
\]

On the other hand, by equation~(2.2.1) and Fact~\ref{2.3}, we have
\[
\begin{aligned}
&\Ext^i_{A\otimes_{\Bbbk}B}
\left(
(A\otimes_{\Bbbk}B)/\mathcal{I},
L\otimes_{\Bbbk}N
\right)
\cong
\bigoplus_{p+q=i}
\Ext^p_A(A/I,L)\otimes_{\Bbbk}
\Ext^q_B(B/J,N).
\end{aligned}
\]
Since the left-hand side vanishes whenever
$i\neq \cd(I,A)+\cd(J,B)$, it follows from $(*)$ and $(**)$ that
$\Ext^s_A(A/I,L)=0$ for all $s>\cd(I,A)$,
and
$\Ext^t_B(B/J,N)=0$ for all $t>\cd(J,B)$.
Therefore,
\[
\Ext^s_A(A/I,L)=0
\quad\text{for all }s\neq \cd(I,A),
\]
and
\[
\Ext^t_B(B/J,N)=0
\quad\text{for all }t\neq \cd(J,B).
\]
Thus, $L$ and $N$ are Gorenstein with respect to $I$ and $J$,
respectively.

$(\mathrm{b})\Rightarrow(\mathrm{a})$:
Suppose that $L$ and $N$ are Gorenstein with respect to $I$ and $J$,
respectively. Then
$\Ext^p_A(A/I,L)=0$ for all $p\neq \cd(I,A)$,
and
$\Ext^q_B(B/J,N)=0$ for all $q\neq \cd(J,B)$.
By equation~(2.2.1) and Fact~\ref{2.3}, it follows that
\[
\Ext^i_{A\otimes_{\Bbbk}B}
\left(
(A\otimes_{\Bbbk}B)/\mathcal{I},
L\otimes_{\Bbbk}N
\right)=0
\quad\text{for all}\quad
i\neq \cd(I,A)+\cd(J,B)
=\cd(\mathcal{I},A\otimes_{\Bbbk}B),
\]
where the last equality follows from Fact~\ref{2.2}(b).
Hence $L\otimes_{\Bbbk}N$ is Gorenstein with respect to
$\mathcal{I}$.
\end{proof}

Let $(A, \fm)$ and $(B, \fn)$ be two local algebras over a field $\Bbbk$, and sharing $\Bbbk$ as the common residue field.
 We also assume that either $A$ or $B$ is algebraic over $\Bbbk$. Then $A\otimes_\Bbbk B$ is a local ring
  with the unique maximal ideal $\mathfrak{M}= \fm \otimes_\Bbbk B + A \otimes_\Bbbk \fn$, see \cite[Remark 2.8]{R4}.
  As a consequence of Theorem \ref{2.4}, we obtain the following corollary.
 \begin{corollary}
\label{2.5}
Under the above assumptions, the following conditions are equivalent:
 \begin{itemize}
\item[{(a)}] $L\otimes_{\Bbbk} N$ is Gorenstein;
\item[{(b)}] Both $L$ and $N$ are Gorenstein.
\end{itemize}
\end{corollary}
\begin{definition}
\label{2.6}
Let $\fa$ be an ideal of $R$, and let $K$ be a finitely generated $R$-module.
We say that $K$ is \emph{maximal Cohen--Macaulay with respect to $\fa$} if $\grade(\fa, K)=\cd(\fa,R)$ (see \cite{PDR}).
\end{definition}
Let $(R,\fm')$ be a local ring. Then $K$ is maximal Cohen--Macaulay if and only if it is maximal Cohen--Macaulay with respect to $\fm'$.
 It is worth noting that if $K$ is Gorenstein with respect to $\fa$ and $\fa K\neq K$, then $K$ is
 maximal Cohen--Macaulay with respect to $\fa$; see also \cite[Theorem~2.19(c)]{PDR}. The converse, however, does not hold in general.
\begin{remark}
\label{2.7}{\em
In view of Fact~\ref{2.2}, the maximal Cohen--Macaulayness of $L\otimes_{\Bbbk} N$ with respect to $\mathcal{I}$ is easily characterized:
 $L\otimes_{\Bbbk} N$ is maximal Cohen--Macaulay with respect to $\mathcal{I}$ if and only if $L$ and $N$ are maximal Cohen--Macaulay with respect to $I$ and $J$, respectively.}
\end{remark}

\section{The arithmetic rank of $L\otimes_{\Bbbk} N$}
 Let $\fa$ be an ideal of $R$ and $K$ a finitely generated $R$-module with $\fa K\neq K$.
 We denote by $\ara(\fa, K)$ the {\it arithmetic rank of $\fa$ with respect to $K$}, which is defined as the infimum of the integers $n \in \mathbb{N}_0$
 such that there exist $x_1,x_2,\ldots,x_n \in R$ satisfying
\[
\Rad\left(\langle x_1,x_2,\ldots,x_n \rangle + \Ann_R K\right) = \Rad\left(\fa+\Ann_R K\right).
\]
It is easy to check that $\cd(\fa, K)\leq \ara(\fa, K)$. Obviously, $\ara(\fa,R)=\ara(\fa)$. If $(R, \fm')$ is a local
ring, then $\ara(\fm')=\ara(\fm', R)=\dim(R)$, and $\ara(\fm', K)=\dim_R(K)$.
For more details and additional information, see \cite{AR2} and \cite{DGTZ}.

In this section, we use the notation and assumptions of the previous section. More precisely, let $\Bbbk$ be a field,
 and let $A$ and $B$ be Noetherian $\Bbbk$-algebras such that $A\otimes_\Bbbk B$ is Noetherian.
  Let $I$ and $J$ be ideals of $A$ and $B$, respectively, and set $\mathcal{I}=I\otimes_\Bbbk B+A\otimes_\Bbbk J$.
   Furthermore, let $L$ and $N$ be nonzero finitely generated modules over $A$ and $B$, respectively.
\begin{fact}	
\label{3.1}{\em Let $\Bbbk$ be an algebraically closed field. Then
\[
\Ass_{A \otimes_\Bbbk B} (L\otimes_\Bbbk N)=\{\fp \otimes_\Bbbk B+ A\otimes_\Bbbk \fq: \fp \in \Ass_A  L \quad\text{and}\quad\fq\in \Ass_B N\},
\]
see \cite[Corollary 2.8]{HNTT}.}
\end{fact}

Recall that the ring $R$ is called \emph{reduced} if it has no nonzero nilpotent elements.
Let $\mathfrak{a}$ be an ideal of a ring $R$.  It is easy to see that the quotient ring $R/\mathfrak{a}$
 is reduced if and only if $\mathfrak{a}$ is a radical ideal, that is,
$
\operatorname{Rad}(\mathfrak{a})=\mathfrak{a}.
$
 A filed $\Bbbk$ is \emph{perfect} if either $\operatorname{char}(\Bbbk)=0$, or $\operatorname{char}(\Bbbk)=p>0$ and every element of $\Bbbk$ has a $p$-th root.
In particular, every algebraically closed field is perfect. Indeed, if $\operatorname{char}(\Bbbk)=0$, then $\Bbbk$ is perfect.
If $\operatorname{char}(\Bbbk)=p>0$, then for every $a\in\Bbbk$, the polynomial $x^p-a$ has a root in $\Bbbk$,
 since $\Bbbk$ is algebraically closed. Thus every element of $\Bbbk$ has a $p$-th root, and hence $\Bbbk$ is perfect.
Finally, if $\Bbbk$ is a perfect field and $A$ and $B$ are reduced $\Bbbk$-algebras, then $A\otimes_{\Bbbk}B$ is reduced, see \cite{B}.

\begin{lemma}\label{3.2}
Let $\Bbbk$ be an algebraically closed field. Then
\[
\Rad(\mathcal{I})
= \Rad(I) \otimes_\Bbbk B + A \otimes_\Bbbk \Rad(J).
\]
\end{lemma}
\begin{proof}
We set
$
H=\Rad(I)\otimes_\Bbbk B+A\otimes_\Bbbk\Rad(J).
$
Let $\fp\in\Ass_A(A/I)$ and $\fq\in\Ass_B(B/J)$. Then
$I\subseteq\fp$ and $J\subseteq\fq$, which imply that
$\Rad(I)\subseteq\fp$ and $\Rad(J)\subseteq\fq$. Consequently,
$
H\subseteq \fp\otimes_\Bbbk B+A\otimes_\Bbbk\fq.
$
Taking intersections over all such $\fp$ and $\fq$, and using
Fact~\ref{3.1}, we obtain
\[
H\subseteq \bigcap_{\substack{\fp \in \Ass_A(A/I) \\ \fq \in \Ass_B(B/J)}}
  \big(\fp \otimes_\Bbbk B + A \otimes_\Bbbk \fq\big)=\Rad(\mathcal I).
\]
On the other hand, since $\mathcal I\subseteq H$, we have
$
\Rad(\mathcal I)\subseteq\Rad(H).
$
Therefore, it remains to show that $H$ is a radical ideal.
Observe that
\[
(A\otimes_\Bbbk B)/H
\cong
(A/\Rad(I))
\otimes_\Bbbk
(B/\Rad(J)).
\]
Since $\Rad(I)$ and $\Rad(J)$ are radical ideals, the quotient
rings $A/\Rad(I)$ and $B/\Rad(J)$ are reduced. Moreover,
$\Bbbk$ is perfect because it is algebraically closed. By the
standard theorem that the tensor product of two reduced algebras
over a perfect field is reduced(see \cite{B}), it follows that
$
(A/\Rad(I))
\otimes_\Bbbk
(B/\Rad(J))
$
is reduced. Hence $(A\otimes_\Bbbk B)/H$ is reduced, which
implies that $H$ is a radical ideal. Thus $\Rad(H)=H$, and
consequently,
$
\Rad(\mathcal I)\subseteq H.
$
Combining the two inclusions, we conclude that
$
\Rad(\mathcal I)
=\Rad(I)\otimes_\Bbbk B
+A\otimes_\Bbbk\Rad(J).
$
This completes the proof.
\end{proof}

\begin{lemma}
\label{3.3}
Let $I\subseteq I'$ be proper ideals of $A$ and
$J\subseteq J'$ be proper ideals of $B$. Suppose that
$
I\otimes_\Bbbk B+A\otimes_\Bbbk J
=
I'\otimes_\Bbbk B+A\otimes_\Bbbk J'.
$
Then $I=I'$ and $J=J'$.
\end{lemma}

\begin{proof}
Consider the canonical surjective map
$
\pi:A\otimes_\Bbbk B
\longrightarrow
(A/I)\otimes_\Bbbk(B/J),
$
defined by
$
\pi(a\otimes_\Bbbk b)
=(a+I)\otimes_\Bbbk(b+J).
$
By the universal property of the tensor product, $\pi$ is
well-defined, and
$
\ker\pi=I\otimes_\Bbbk B+A\otimes_\Bbbk J.
$
By the assumed equality,
$
I'\otimes_\Bbbk B+A\otimes_\Bbbk J'
=\ker\pi.
$
Hence
$I'\otimes_\Bbbk B\subseteq\ker\pi$ and $A\otimes_\Bbbk J'\subseteq\ker\pi$.
Consequently,
$
(I'/I)\otimes_\Bbbk(B/J)=0$ and $(A/I)\otimes_\Bbbk(J'/J)=0$.
Since $I$ and $J$ are proper ideals, the $\Bbbk$-vector spaces
$A/I$ and $B/J$ are nonzero. As the tensor product of two
nonzero vector spaces over a field is nonzero, we obtain
$I'/I=0$ and  $J'/J=0$.
Therefore, $I=I'$ and $J=J'$.
\end{proof}

\begin{theorem}
\label{3.4}
Let $\Bbbk$ be an algebraically closed field. Then
\[
\ara(\mathcal{I}, L \otimes_\Bbbk N)
= \ara(I, L) + \ara(J, N).
\]
\end{theorem}

\begin{proof}
Let $\ara(I, L) = r$ and $\ara(J, N) = s$. Then there exist
$x_1, \dots, x_r \in A$ and $y_1, \dots, y_s \in B$ such that
\[
\Rad(\langle x_1, \dots, x_r \rangle + \Ann_A L)= \Rad(I + \Ann_A L),
\]
and
\[
\Rad(\langle y_1, \dots, y_s \rangle + \Ann_B N)= \Rad(J + \Ann_B N).
\]
Set $x =  x_1, \dots, x_r $ and  $y = y_1, \dots, y_s$.
We claim that
\[
\Rad\!\big(
\langle x \rangle \otimes_\Bbbk B
+ A \otimes_\Bbbk \langle y \rangle
+ \Ann_{A \otimes_\Bbbk B}(L \otimes_\Bbbk N)
\big)
= \Rad\!\big(
\mathcal{I} + \Ann_{A \otimes_\Bbbk B}(L \otimes_\Bbbk N)
\big).
\]
Note that
\[
\Ann_{A \otimes_\Bbbk B}(L \otimes_\Bbbk N)
= (\Ann_A L) \otimes_\Bbbk B
+ A \otimes_\Bbbk \Ann_B N\  \  \   (*).
\]
We set $\Ann_{A \otimes_\Bbbk B}(L \otimes_\Bbbk N)=K$.
Observe that
\begin{gather*}
\begin{aligned}
\Rad\!\big(
\langle x \rangle \otimes_\Bbbk B
+ A \otimes_\Bbbk \langle y \rangle
+K
\big)
&= \Rad\!\big(
(\langle x \rangle + \Ann_A L) \otimes_\Bbbk B
+ A \otimes_\Bbbk (\langle y \rangle + \Ann_B N)
\big) \\[2pt]
&= \Rad(\langle x \rangle + \Ann_A L) \otimes_\Bbbk B
+ A \otimes_\Bbbk \Rad(\langle y \rangle + \Ann_B N) \\[2pt]
&= \Rad(I + \Ann_A L) \otimes_\Bbbk B
+ A \otimes_\Bbbk \Rad(J + \Ann_B N) \\[2pt]
&= \Rad\!\big(
I \otimes_\Bbbk B
+ (\Ann_A L) \otimes_\Bbbk B
+ A \otimes_\Bbbk J
+ A \otimes_\Bbbk \Ann_B N
\big) \\[2pt]
&= \Rad\!\big(
\mathcal{I}
+ \Ann_{A \otimes_\Bbbk B}(L \otimes_\Bbbk N)
\big).
\end{aligned}
\end{gather*}

Consequently,
\[
\ara(\mathcal{I}, L \otimes_\Bbbk N) \le r + s.
\]
Here, the second equality follows from Lemma \ref{3.2}, while the remaining steps are standard.

Now let $\ara(\mathcal{I}, L \otimes_\Bbbk N) = t$.
Then there exist elements $x_1, \dots, x_r \in A$
and $y_1, \dots, y_s \in B$ with $r + s = t$ such that
\[
\Rad\!\big(
\langle x \rangle \otimes_\Bbbk B
+ A \otimes_\Bbbk \langle y \rangle
+ \Ann_{A \otimes_\Bbbk B}(L \otimes_\Bbbk N)
\big)
= \Rad\!\big(
\mathcal{I} + \Ann_{A \otimes_\Bbbk B}(L \otimes_\Bbbk N)
\big).
\]
Hence, in view of  $(*)$ and Lemma \ref{3.2}, we have
\[
\Rad(\langle x \rangle + \Ann_A L) \otimes_\Bbbk B
+ A \otimes_\Bbbk \Rad(\langle y \rangle + \Ann_B N)
= \Rad(I + \Ann_A L) \otimes_\Bbbk B
+ A \otimes_\Bbbk \Rad(J + \Ann_B N).
\]
It follows from Lemma \ref{3.3} that
\[
\Rad(\langle x \rangle + \Ann_A L)= \Rad(I + \Ann_A L),
\]
and
\[
\Rad(\langle y \rangle + \Ann_B N)= \Rad(J + \Ann_B N).
\]
Thus $\ara(I, L) \le r$ and $\ara(J, N) \le s$, and therefore
\[
\ara(I, L) + \ara(J, N)
\le r + s = t
= \ara(\mathcal{I}, L \otimes_\Bbbk N),
\]
as desired.
\end{proof}

\section{Approximate Cohen-Macaulayness of $L\otimes_{\Bbbk} N$}
In this section, let $\fa$ an ideal of $R$ and let $K$ be a finitely generated $R$-module.
Moreover, we keep the notation and assumptions introduced in the previous sections.
\begin{definition}
\label{4.1} We say $K$ is {\it Cohen--Macaulay with respect to $\fa$} if either $K=\fa K$ or $K\neq \fa K$
and $\grade(\fa, K)=\cd(\fa, K)$ (see \cite{R1}).
\end{definition}
We note that if $(R, \fm')$ is a local ring,
 then $K$ is Cohen–Macaulay if and only if $K$ is Cohen–Macaulay with respect to $\fm'$.

\begin{definition}
\label{4.2}
We say $K$ is \emph{approximately Cohen-Macaulay module with respect to \( \fa \)} if either $\cd(\fa, K)=0$ or there exists an element $a\in \fa$
 such that $K/a^nK$ is Cohen-Macaulay with respect to $\fa $ with $\cd(\fa, K/a^nK)=\cd(\fa, K)-1$, for every integer $n\geq 1$(see \cite{FR}).
\end{definition}
Note that if \((R,\fm')\) is a local ring, then \(K\) is approximately Cohen--Macaulay if and only if \(K\) is approximately Cohen--Macaulay with respect to \(\fm'\); see \cite{G,S}.
One observes that if \(K\) is Cohen--Macaulay with respect to \(\fa\), then \(K\) is approximately Cohen--Macaulay with respect to \(\fa\); see \cite[Proposition 2.4]{FR}.

\begin{definition}
\label{4.3} {\em A filtration $\mathscr{D}$: $0=D_0\varsubsetneq D_1\varsubsetneq \dots  \varsubsetneq  D_r=K$ of $K$ by submodules $D_i$ is called the
 {\em dimension filtration of $K$ with respect to $\fa$} if $D_{i-1}$ is the largest submodule of $D_i$ for which $\cd(\fa, D_{i-1})<\cd(\fa, D_i)$
for all $i=1, \dots, r$.}
\end{definition}
A finite filtration $\mathcal{F}: 0=K_0\varsubsetneq K_1 \varsubsetneq \dots  \varsubsetneq K_r=K$ of $K$ by submodules $K_i$ is called a
{\em Cohen-Macaulay filtration with respect to $\fa$} if each quotient $K_i/K_{i-1}$ is Cohen-Macaulay with respect to $\fa$ and
$0 \leq \cd(\fa, K_1/K_0)<\cd(\fa, K_2/K_1)< \dots< \cd(\fa, K_r/K_{r-1})$. If $K$ admits a Cohen-Macaulay filtration with respect to $\fa$,
then we say $K$ is \emph {sequentially Cohen-Macaulay with respect to $\fa$}. Note that if $K$ is sequentially Cohen-Macaulay with respect to $\fa$,
then the filtration $\mathscr{F}$ is uniquely determined and it is just the dimension filtration of $K$ with respect to $\fa$,
 i.e., $\mathscr{F}=\mathscr{D}$, see Section 4 of \cite{R4}.
\begin{fact}
\label{4.4}{\em
Let $\fa$ be an ideal of $R$ contained in the Jacobson radical of $R$ with $\cd(\fa, K) = \ara(\fa, K)$.
 Suppose that $K$ is not Cohen–Macaulay with respect to $\fa$. Then the following statements are equivalent:
\begin{enumerate}
\item [(a)] $K$ is approximately Cohen–Macaulay with respect to $\fa$;
\item [(b)] $K$ is sequentially Cohen–Macaulay with respect to $\fa$ with the dimension filtration $\mathcal{F}$:
$0 = D_0 \varsubsetneq D_1 \varsubsetneq D_2 = K$ with respect to $\fa$,
where $\cd(\fa, D_1) = \cd(\fa, K) - 1$, see \cite[Theorem 3.2]{FR}.
\end{enumerate}
%Here, by saying that $M$ is sequentially Cohen--Macaulay with respect to $\fa$, we mean that $D_1/D_0$ and $D_2/D_1$ are Cohen--Macaulay with respect to $\fa$.
In particular, let $(R, \fm')$ be a local ring, and $K$ a non-Cohen--Macaulay $R$-module.
  Then, $K$ is approximately Cohen-Macaulay if and only if
$K$ is sequentially Cohen-Macaulay  with the dimension filtration $\mathcal{F}$:
 $0=D_0\varsubsetneq D_1 \varsubsetneq D_2=K$ where $\dim D_1=\dim_RK-1$.
}
\end{fact}

\begin{example}
\label{4.5}
{\em
Recall that $L$ is a sequentially Cohen--Macaulay $A$-module and $N$ is a sequentially Cohen--Macaulay $B$-module
if and only if $L \otimes_\Bbbk N$ is a sequentially Cohen--Macaulay $A\otimes_\Bbbk B$-module, see \cite[Theorem 4.6]{R4}.
However, this equivalence does not hold for approximately Cohen--Macaulay modules.
Indeed, set
\[
L=\Bbbk[x_1, x_2, x_3]/\big(\langle x_1 \rangle \cap \langle x_2, x_3 \rangle\big),
\]
viewed as a $\Bbbk[x]=\Bbbk[x_1, x_2, x_3]$-module.
One verifies that $\dim L = 2$ and $\depth L = 1$.
We claim that $L$ is an approximately Cohen--Macaulay $\Bbbk[x]$-module.
Consider the dimension filtration $0 \subsetneq L_1 \subsetneq L$, where
\[
L_1 = \langle x_1 \rangle / \big(\langle x_1 \rangle \cap \langle x_2, x_3 \rangle\big).
\]
Then $L/L_1 \cong \Bbbk[x]/\langle x_1 \rangle$
is Cohen--Macaulay of dimension $2$. Moreover,
\[
\Ass_{\Bbbk[x]}(L_1) = \{\langle x_2, x_3 \rangle\},
\]
so $L_1$ is Cohen--Macaulay of dimension $1$.
Similarly,
\[
N=\Bbbk[y_1, y_2, y_3]/\big(\langle y_1 \rangle \cap \langle y_2, y_3 \rangle\big)
\]
is an approximately Cohen--Macaulay $\Bbbk[y]$-module.
On the other hand, $L \otimes_\Bbbk N$ is \emph{not} approximately Cohen--Macaulay as an $\Bbbk[x, y]$-module, because
\[
\dim (L \otimes_\Bbbk N) = \dim L + \dim N = 2+2=4,
\]
while
\[
\depth (L \otimes_\Bbbk N) = \depth L + \depth N = 1+1=2.
\]
Hence, $\dim (L \otimes_\Bbbk N) - \depth (L \otimes_\Bbbk N) = 2$,
which contradicts approximate Cohen--Macaulayness.
}
\end{example}
\begin{remark}{\em
\label{4.6}
%Suppose $IL\neq L$ and $JN\neq N$.
 In view of Fact \ref{2.2}, $L\otimes_{\Bbbk} N$ is Cohen--Macaulay with respect to $\mathcal{I}$
 if and only if $L$ and $N$ are Cohen--Macaulay with respect to $I$ and $J$, respectively.}
\end{remark}

\begin{theorem}
\label{4.7}
Let $\Bbbk$ be an algebraically closed field. Suppose $\cd(I, L)=\ara(I, L)$ and $\cd(J, N)=\ara(J, N)$.
 The following statements are equivalent:
\begin{itemize}
\item[(a)] $L \otimes_\Bbbk N$ is approximately Cohen--Macaulay with respect to $\mathcal{I}$ which is not Cohen-Macaulay with respect to $\mathcal{I}$.
\item[(b)] $L$ is approximately Cohen--Macaulay with respect to $I$ which is not Cohen--Macaulay with respect to $I$, and $N$ is Cohen--Macaulay with respect to $J$.
\end{itemize}
\end{theorem}
\begin{prf}
\((a) \Rightarrow (b)\): Suppose $L\otimes_\Bbbk N$ is approximately Cohen--Macaulay with respect to $\mathcal{I}$ but not Cohen-Macaulay with respect to $\mathcal{I}$.
 Let  $\mathscr{F}$: $0=L_0\varsubsetneq L_1 \varsubsetneq \dots  \varsubsetneq L_r=L$
 be the dimension filtration of $L$ with respect to $I$. A proof similar to that of \cite[Lemma 4.5]{R4} shows that
\[
\mathscr{F}\otimes_{\Bbbk}N:\quad
0=L_0\otimes_{\Bbbk}N
\varsubsetneq L_1\otimes_{\Bbbk}N
\varsubsetneq \dots
\varsubsetneq L_r\otimes_{\Bbbk}N
=L\otimes_{\Bbbk}N
\]
is the dimension filtration of \(L\otimes_{\Bbbk}N\) with respect to \(\mathcal{I}\).
By our assumption $L\otimes_\Bbbk N$ is sequentially Cohen--Macaulay with respect to $\mathcal{I}$ with the Cohen--Macaulay filtration of length $2$, see Fact \ref{4.4}.
 Note that the Cohen--Macaulay filtration is the same as dimension filtration. Thus, $0 \subsetneq L_1\otimes_\Bbbk N \subsetneq L\otimes_\Bbbk N$ is
  the Cohen--Macaulay filtration of $L\otimes_\Bbbk N$ with respect to $\mathcal{I}$ with $\cd(\mathcal{I}, L_1\otimes_\Bbbk N)=\cd(\mathcal{I}, L\otimes_\Bbbk N)-1$.
   This yields $0 \subsetneq L_1 \subsetneq L$ is a Cohen--Macaulay filtration for $L$ with $\cd(I,  L_1)=\cd(I, L)-1$, and $N$ is Cohen--Macaulay with respect to $J$. In fact, as
   $L_1\otimes_\Bbbk N$ is Cohen--Macaulay with respect $\mathcal{I}$, it follows from Remark \ref{4.6} that $L_1$ is Cohen--Macaulay with respect to $I$, and $N$ is Cohen--Macaulay with respect to $J$. Moreover, as
   \[
(L\otimes_\Bbbk N)/(L_1\otimes_\Bbbk N)\cong L/L_1\otimes_\Bbbk N
\]
is Cohen--Macaulay with respect $\mathcal{I}$, it follows again from Remark \ref{4.6} that $L_2/L_1$ is Cohen--Macaulay with respect to $I$.
Finally, by Fact \ref{2.2}(b), the equality of $\cd(\mathcal{I}, L_1\otimes_\Bbbk N)=\cd(\mathcal{I}, L\otimes_\Bbbk N)-1$ implies that $\cd(I,  L_1)=\cd(I, L)-1$.
 Therefore, $L$ is approximately Cohen--Macaulay with respect to $I$ by Fact \ref{4.4}, and $N$ is Cohen--Macaulay with respect to $J$.
Note also that $L$ is not Cohen--Macaulay with respect to $I$. Indeed, if $L$ were Cohen--Macaulay with respect to $I$,
 then $\grade(I,L)=\cd(I,L)$. Since $N$ is Cohen--Macaulay with respect to $J$, Fact~\ref{2.2} would then imply that
$
\grade(\mathcal{I},L\otimes_{\Bbbk}N)
=\cd(\mathcal{I},L\otimes_{\Bbbk}N),
$
yielding a contradiction.

\((b) \Rightarrow (a)\): Assume that $L$ is approximately Cohen--Macaulay with respect to $I$ but not Cohen--Macaulay with respect to $I$, and $N$ is Cohen--Macaulay with respect to $J$.
  Thus, there exists Cohen--Macaulay filtration of $0  \subsetneq L_1 \subsetneq L$ with respect to $I$ with
$\cd(I, L_1)=\cd(I,  L)-1$.
Our assumption, Fact \ref{2.2} and Theorem \ref{3.4} yield
\[
\cd(\mathcal{I}, L\otimes_\Bbbk N)= \cd(I, L)+\cd(J, N)= \ara(I, L)+\ara(J, N)=\ara(\mathcal{I}, L\otimes_\Bbbk N).
\]
We claim that $0\subsetneq L_1\otimes_\Bbbk N \subsetneq L\otimes_\Bbbk N$ is a Cohen--Macaulay filtration for $L\otimes_\Bbbk N$ with
$\cd(\mathcal{I}, L_1\otimes_\Bbbk N) =\cd(\mathcal{I}, L\otimes_\Bbbk N)-1$.
In fact, as $L_1$ is Cohen--Macaulay with respect to $I$ and $N$ is Cohen--Macaulay with respect to $J$, it follows from Remark \ref{4.6} that $L_1\otimes_\Bbbk N$ is Cohen--Macaulay with respect to $\mathcal{I}$. Moreover,
\[
(L\otimes_\Bbbk N)/(L_1\otimes_\Bbbk N)\cong L/L_1\otimes_\Bbbk N
\]
is Cohen--Macaulay with respect to $\mathcal{I}$ again by Remark \ref{4.6}. Finally, the equality of $\cd(I,  L_1)=\cd(I, L)-1$
  implies that $\cd(\mathcal{I}, L_1\otimes_\Bbbk N)=\cd(\mathcal{I}, L\otimes_\Bbbk N)-1$.
 Therefore,  $L \otimes_\Bbbk N$ is approximately Cohen--Macaulay with respect to $\mathcal{I}$ by Fact \ref{4.4}.
Note also that $L\otimes_\Bbbk N$ is not Cohen--Macaulay with respect to $\mathcal{I}$. Otherwise, by Fact~\ref{2.2}, $L$ would be Cohen--Macaulay with respect to $I$, a contradiction.
 \end{prf}

 Let $(A, \fm)$ and $(B, \fn)$ be two local algebras over a field $\Bbbk$, and sharing $\Bbbk$ as the common residue field.
 We also assume that either $A$ or $B$ is algebraic over $\Bbbk$. Then $A\otimes_\Bbbk B$ is a local ring
  with the unique maximal ideal $\mathfrak{M}= \fm \otimes_\Bbbk B + A \otimes_\Bbbk \fn$, see \cite{R4}.
  Note that $\cd(\fm, L)=\ara(\fm, L)=\dim_AL$ and $\cd(\fn, N)=\ara(\fn, N)=\dim_BN$.
    Consequently, Theorem \ref{4.7} yields the following corollary.
 \begin{corollary}
\label{4.8}
Let $\Bbbk$ be an algebraically closed field. Under the above assumptions, the following conditions are equivalent:
 \begin{itemize}
\item[{(a)}] $L\otimes_\Bbbk N$ is approximately Cohen--Macaulay which is not Cohen--Macaulay.
\item[{(b)}] $L$ is approximately Cohen--Macaulay and $N$ is Cohen--Macaulay.
 \end{itemize}
 \end{corollary}
Setting $I=0$ in Theorem \ref{4.7}, we obtain the following corollary.
\begin{corollary}
\label{4.9}
Let $\Bbbk$ be an algebraically closed field, and suppose that
$\cd(J,N)=\ara(J,N)$. The following conditions are equivalent:
\begin{itemize}
    \item[(a)] $L\otimes_\Bbbk N$ is approximately Cohen--Macaulay but not Cohen--Macaulay with respect to $A\otimes_\Bbbk J$.
    \item[(b)] $N$ is approximately Cohen--Macaulay but not Cohen--Macaulay with respect to $J$.
\end{itemize}
\end{corollary}

\section{Maximal depth property of $L\otimes_{\Bbbk} N$}
Let $\fa$ be an ideal of $R$ and $K$ a finitely generated $R$-module.
Observe that
\[
\grade(\fa, K)\leq \min\{\cd(\fa, R/\fp): \fp\in \Ass_RK\},
\]
see \cite[Proposition 2.2]{R6}. Motivated by this inequality, we define
\[
\mcd(\fa, K)=\min\{\cd(\fa, R/\fp): \fp\in \Ass_RK\}.
\]
Hence, $\grade(\fa, K)\leq \mcd(\fa, K)$.
We say that $K$ has \emph{maximal depth with respect to $\fa$} if $\grade(\fa, K)=\mcd(\fa, K)$.
Equivalently, $K$ has maximal depth with respect to $\fa$ if there exists
$\fp\in\Ass_RK$ such that $\grade(\fa, K)=\cd(\fa, R/\fp)$. If $(R,\fm')$ is a local ring, then
\[
\grade(\fm', K)=\depth_R K
\quad\text{and}\quad
\mcd(\fm', K)=\mcd_R K.
\]
Therefore, $K$ has \emph{maximal depth} if and only if it has maximal depth with respect to $\fm'$, that is,
$\depth_R K=\mcd_R K$.
Equivalently, $K$ has maximal depth if and only if there exists
$\fp\in\Ass_RK$ such that $\depth_R K=\dim R/\fp$. Note that
\[
\grade(\fa, K)=0
\quad\Longleftrightarrow\quad
\mcd(\fa, K)=0.
\]
Moreover, if $\fa K\neq K$, then
 \begin{eqnarray}	
\label{5.0.1}
\grade(\fa, K)\leq \mcd(\fa, K)\leq \max \{\cd(\fa, R/\fp): \fp\in \Ass_RK\}=\cd(\fa, K).
\end{eqnarray}

The last equality follows from \cite[Corollary~4.6]{CJR}.
A proof similar to that of \cite[Proposition 1.6]{R3} shows that if $K$ is sequentially Cohen--Macaulay with respect to $\fa$,
 then $K$ has maximal depth with respect to $\fa$. In particular, by Fact \ref{4.4}, if $K$ is approximately Cohen--Macaulay
  with respect to $\fa$, then $K$ has maximal depth with respect to $\fa$.
 We keep the notation and assumptions introduced in the previous sections.
\begin{proposition}
\label{5.1} Let $\Bbbk$ be an algebraically closed field. Then
\[
\mcd(\mathcal{I}, L\otimes_\Bbbk N)
=
\mcd(I,L)+\mcd(J,N).
\]
\end{proposition}
\begin{prf}
Let $\mcd(I,L)=s$ and $\mcd(J,N)=t$. Then there exist
$\fp_1\in\Ass_A L$ and $\fp_2\in\Ass_B N$ such that
\[
\cd(I,A/\fp_1)=s
\quad\text{and}\quad
\cd(J,B/\fp_2)=t.
\]
By Fact~\ref{3.1},
$
\fp=\fp_1\otimes_\Bbbk B+A\otimes_\Bbbk\fp_2
$
belongs to $\Ass_{A\otimes_\Bbbk B}(L\otimes_\Bbbk N)$. Moreover,
\[
\cd(\mathcal{I},(A\otimes_\Bbbk B)/\fp)
=
\cd(I,A/\fp_1)+\cd(J,B/\fp_2)
=s+t.
\]
Therefore,
$
\mcd(\mathcal{I},L\otimes_\Bbbk N)\leq s+t.
$

Conversely, let
$
\mcd(\mathcal{I},L\otimes_\Bbbk N)=k.
$
Then there exists
$\fp\in\Ass_{A\otimes_\Bbbk B}(L\otimes_\Bbbk N)$ such that
$
\cd(\mathcal{I},(A\otimes_\Bbbk B)/\fp)=k.
$
By Fact~\ref{3.1}, $\fp=\fp_1\otimes_\Bbbk B+A\otimes_\Bbbk\fp_2$ for some $\fp_1\in\Ass_A L$ and $\fp_2\in\Ass_B N$. Hence,
\[
\cd(I,A/\fp_1)+\cd(J,B/\fp_2)=k.
\]
By the definition of $\mcd$, we have
\[
\mcd(I,L)\leq\cd(I,A/\fp_1)
\quad\text{and}\quad
\mcd(J,N)\leq\cd(J,B/\fp_2).
\]
Consequently,
\[
\mcd(I,L)+\mcd(J,N)
\leq
\cd(I,A/\fp_1)+\cd(J,B/\fp_2)
=k.
\]
Thus,
\[
\mcd(I,L)+\mcd(J,N)
\leq
\mcd(\mathcal{I},L\otimes_\Bbbk N).
\]
Combining the two inequalities gives the desired equality.
\end{prf}

\begin{theorem}
\label{5.2}
Let $\Bbbk$ be an algebraically closed field. Then the following conditions are equivalent:
\begin{itemize}
\item[(a)] $L\otimes_{\Bbbk} N$ has maximal depth with respect to $\mathcal{I}$;
\item[(b)] $L$ and $N$ have maximal depth with respect to $I$ and $J$, respectively.
\end{itemize}
\end{theorem}

\begin{prf}
$(\mathrm{a})\Rightarrow(\mathrm{b})$: Suppose that $L\otimes_{\Bbbk} N$ has maximal depth with respect to $\mathcal{I}$. Then
\[
\grade(\mathcal{I},L\otimes_\Bbbk N)
=
\mcd(\mathcal{I},L\otimes_\Bbbk N).
\]
Since
$
\grade(\mathcal{I},L\otimes_\Bbbk N)
=
\grade(I,L)+\grade(J,N)
$
and, by Proposition~\ref{5.1},
\[
\mcd(\mathcal{I},L\otimes_\Bbbk N)
=
\mcd(I,L)+\mcd(J,N),
\]
we obtain
\[
\grade(I,L)+\grade(J,N)
=
\mcd(I,L)+\mcd(J,N).
\]
Since
\[
\grade(I,L)\leq\mcd(I,L)
\quad\text{and}\quad
\grade(J,N)\leq\mcd(J,N),
\]
it follows that
\[
\grade(I,L)=\mcd(I,L)
\quad\text{and}\quad
\grade(J,N)=\mcd(J,N).
\]
Hence, $L$ and $N$ have maximal depth with respect to $I$ and $J$, respectively.

$(\mathrm{b})\Rightarrow(\mathrm{a})$: Suppose that $L$ and $N$ have maximal depth with respect to $I$ and $J$, respectively. Then
\[
\grade(I,L)=\mcd(I,L)
\quad\text{and}\quad
\grade(J,N)=\mcd(J,N).
\]
Therefore,
\[
\begin{aligned}
\grade(\mathcal{I},L\otimes_\Bbbk N)
&=\grade(I,L)+\grade(J,N)\\
&=\mcd(I,L)+\mcd(J,N)\\
&=\mcd(\mathcal{I},L\otimes_\Bbbk N),
\end{aligned}
\]
where the last equality follows from Proposition~\ref{5.1}. Thus,
$L\otimes_\Bbbk N$ has maximal depth with respect to $\mathcal{I}$.
\end{prf}

Let $(A, \fm)$ and $(B, \fn)$ be two Noetherian local algebras over a field $\Bbbk$ such that $A \otimes_\Bbbk B$ is Noetherian.  We also suppose
that $\Bbbk$ is the common residue field of both $A$ and $B$, and either $A$ or $B$ is algebraic over $\Bbbk$.
As a consequence of Theorem \ref{5.2}, we obtain the following result.
\begin{corollary}
\label{5.3}
Under the above assumptions, suppose further that $\Bbbk$ is algebraically closed. The following conditions are equivalent:
\begin{itemize}
\item[(a)] $L\otimes_{\Bbbk} N$ has maximal depth;
\item[(b)] $L$ and $N$ have maximal depth.
\end{itemize}
\end{corollary}

\begin{prf}
Set
$
\mathfrak{M}=\mathfrak{m}\otimes_{\Bbbk} B
+A\otimes_{\Bbbk}\mathfrak{n}.
$
Then \(\mathfrak{M}\) is the unique maximal ideal of
\(A\otimes_{\Bbbk}B\).
Part (a) means that \(L\otimes_{\Bbbk}N\) has maximal depth
with respect to \(\mathfrak{M}\), while part (b) means that
\(L\) and \(N\) have maximal depth with respect to
\(\mathfrak{m}\) and \(\mathfrak{n}\), respectively.
Therefore, the assertion follows immediately from
Theorem~\ref{5.2}.
\end{prf}

\section{Homological properties of $M/(I\otimes_{\Bbbk}B)M$}
In this section, we use the notation and assumptions introduced in the previous sections.
In the following, we investigate similar questions raised in \cite[Question 6.4]{R5} (see also \cite[Question 1.2]{NR})
concerning $L\otimes_\Bbbk N$. In fact, we show that if $M=L\otimes_\Bbbk N$ is Cohen--Macaulay,
generalized Cohen--Macaulay, approximately Cohen--Macaulay, or sequentially Cohen--Macaulay with
respect to $A\otimes_\Bbbk J$, then $M/(I\otimes_\Bbbk B)M$ has the corresponding property with
respect to $\mathcal{I}$. Moreover, if $M$ has maximal depth with respect to
$A\otimes_\Bbbk J$, then $M/(I\otimes_\Bbbk B)M$ has maximal depth with respect to $\mathcal{I}$. In contrast, Gorensteinness
and maximal Cohen--Macaulayness are not preserved in general under this quotient construction.

As an application, we recover Propositions 6.1--6.3 of \cite{R5}.
Indeed, let
$
S=\Bbbk[x_1,\ldots,x_m,y_1,\ldots,y_n]
$
be the standard bigraded polynomial ring over the field $\Bbbk$,
and set $P=\langle x_1, \ldots, x_m \rangle$ and $Q= \langle y_1, \ldots, y_n \rangle$.
Let $A=\Bbbk[x_1, \ldots, x_m]$ and $B=\Bbbk[y_1,\ldots,y_n]$,
so that $S=A\otimes_\Bbbk B$.
We show that an $S$-module $M$ is Cohen--Macaulay
(respectively, generalized Cohen--Macaulay or sequentially
Cohen--Macaulay) with respect to $Q$ if and only if $M/PM$
has the corresponding property with respect to $P+Q$.

We begin by considering the Cohen--Macaulayness.

\vspace{.2cm}
\noindent{\bf Cohen--Macaulayness:}
\vspace{.3cm}
\begin{fact}
\label{6.1}{\em
By setting $I=0$ in Fact~\ref{2.2} and using
$\grade(0, K)=\cd(0, K)=0$ for every nonzero $R$-module $K$,
we obtain the following equalities.
\begin{itemize}
\item[(a)] $\grade(A\otimes_\Bbbk J, L\otimes_\Bbbk N)=\grade(J,N).$
\item[(b)] $\cd(A\otimes_\Bbbk J, L\otimes_\Bbbk N)=\cd(J,N).$
\end{itemize}}
\end{fact}
\begin{proposition}
\label{6.2}
Set $M = L \otimes_\Bbbk N$.  The following statements are equivalent:
\begin{itemize}
\item[(a)] $M$ is Cohen--Macaulay with respect to $A \otimes_\Bbbk J$;
\item[(b)] $M /(I \otimes_\Bbbk B)M$ is Cohen--Macaulay with respect to $\mathcal{I}$.
\end{itemize}
\end{proposition}
\begin{prf}
Note that
$$\begin{array}{ll}
M/(I \otimes_\Bbbk B)M &\cong ((A \otimes_\Bbbk B)/(I \otimes_\Bbbk B))\otimes_{A \otimes_\Bbbk B} (L\otimes_\Bbbk N)\\
                       & \cong  (A/I\otimes_\Bbbk B/0)\otimes_{A \otimes_\Bbbk B} (L\otimes_\Bbbk N)\\
                       & \cong (A/I\otimes_A L)\otimes_\Bbbk (B\otimes_BN)\\
                       & \cong (L/IL) \otimes_\Bbbk N.
\end{array}$$
Consequently,
 \begin{eqnarray}	
\label{6.2.1}
M/(I \otimes_\Bbbk B)M \cong (L/IL) \otimes_\Bbbk N.
\end{eqnarray}
 Since $L/IL$ is $I$-torsion, it follows that $\cd(I, L/IL)=0$. Thus, by Fact \ref{2.2} and Fact \ref{6.1}, we have
$$\begin{array}{ll}
\grade(\mathcal{I}, M/(I \otimes_\Bbbk B)M) &=\grade(\mathcal{I}, (L/IL)\otimes_\Bbbk N)\\
                                            &= \grade(I, L/IL)+\grade(J, N)\\
                                            & =\grade(J, N)\\
                                            &=\grade(A \otimes_\Bbbk J, M)
\end{array}$$
and
$$\begin{array}{ll}
\cd(\mathcal{I}, M/(I \otimes_\Bbbk B)M)& =\cd(\mathcal{I}, (L/IL)\otimes_\Bbbk N)\\
                                        &= \cd(I, L/IL)+\cd(J, N)\\
                                        &=\cd(J, N)\\
                                        &=\cd(A \otimes_\Bbbk J, M).
\end{array}$$
 Therefore, the assertions  follow.
\end{prf}

\vspace{.2cm}
\noindent{\bf Generalized Cohen--Macaulayness:}
\vspace{.3cm}

\begin{fact}{\em
\label{6.3}
By \cite[Theorem 3.2]{R4}, we have
\[
\f_{\mathcal{I}}(L\otimes_{\Bbbk} N)
=
\min\bigl\{
\grade(I,L)+\f_J(N),\,
\f_I(L)+\grade(J,N)
\bigr\}.
\]
In particular, setting \(I=0\), we have
\(\grade(I, L)=0\) and \(\f_I(L)=+\infty\). Hence,
$
\f_{A\otimes_{\Bbbk} J}(L\otimes_{\Bbbk} N)
=\f_J(N).
$
}
\end{fact}

\begin{proposition}
\label{6.4}
Set $M=L\otimes_\Bbbk N$. The following statements are equivalent:
\begin{itemize}
\item[(a)] $M$ is generalized Cohen--Macaulay with respect to
$A\otimes_\Bbbk J$;
\item[(b)] $M/(I\otimes_\Bbbk B)M$ is generalized Cohen--Macaulay
with respect to $\mathcal{I}$.
\end{itemize}
\end{proposition}

\begin{proof}
Observe that
\begin{align*}
\f_{\mathcal{I}}\big(M/(I\otimes_\Bbbk B)M\big)
&=\f_{\mathcal{I}}\big((L/IL)\otimes_\Bbbk N\big)\\
&=\min\big\{
\grade(I,L/IL)+\f_J(N),
\f_I(L/IL)+\grade(J,N)
\big\}\\
&=\f_J(N)\\
&=\f_{A\otimes_\Bbbk J}(L\otimes_\Bbbk N)=\f_{A\otimes_\Bbbk J}(M).
\end{align*}
The first equality follows from equation (\ref{6.2.1}). The remaining equalities follow from Fact \ref{6.3}. Note also that $\grade(I,L/IL)=0$ and $\f_I(L/IL)=+\infty$.

Moreover, by Fact \ref{2.2}(b) and Fact \ref{6.1}, we have
\[
\begin{aligned}
\cd\big(\mathcal{I}, M/(I\otimes_\Bbbk B)M\big)
&=\cd\big(\mathcal{I},(L/IL)\otimes_\Bbbk N\big)\\
&=\cd(J,N)\\
&=\cd(A\otimes_\Bbbk J, L\otimes_\Bbbk N)=\cd(A\otimes_\Bbbk J, M).
\end{aligned}
\]
 Hence, the assertion follows.
\end{proof}

\vspace{.2cm}
\noindent{\bf Sequential Cohen--Macaulayness:}
\vspace{.3cm}
\begin{fact}{\em
\label{6.5}
Let $\Bbbk$ be an algebraically closed field. Then the following conditions are equivalent:
 \begin{itemize}
\item[{(a)}] $L\otimes_{\Bbbk} N$ is sequentially Cohen--Macaulay with respect to $\mathcal{I}$;
\item[{(b)}] $L$ and $N$ are sequentially Cohen--Macaulay with respect to $I$ and $J$, respectively.
\end{itemize}
See \cite[Theorem 4.6]{R4}.
In particular, we have  $L\otimes_\Bbbk N$ is sequentially Cohen--Macaulay with respect
to $A \otimes_\Bbbk J$ if and only if $N$ is sequentially Cohen--Macaulay with respect to $J$,
see \cite[Corollary 4.7]{R4}.
 }
\end{fact}

\begin{proposition}
\label{6.6}
We set $M=L\otimes_\Bbbk N$. Let $\Bbbk$ be an algebraically closed
field. Then the following statements are equivalent:
 \begin{itemize}
\item[{(a)}]  $M$ is a sequentially Cohen--Macaulay with respect to $A \otimes_\Bbbk J$;
\item[{(b)}] $M/(I \otimes_\Bbbk B)M$ is sequentially Cohen--Macaulay with respect to $\mathcal{I}$.
 \end{itemize}
\end{proposition}
\begin{prf}
By equation (\ref{6.2.1}), we have $M/(I \otimes_\Bbbk B)M\cong L/IL\otimes_\Bbbk N$.
 Since $\cd(I, L/IL)=0$, the module $L/IL$ is Cohen--Macaulay with respect to $I$ and hence sequentially Cohen--Macaulay with respect to $I$.
 Therefore, the desired equivalence follows from Fact \ref{6.5}.
\end{prf}
\vspace{.2cm}
\noindent{\bf Approximate Cohen--Macaulayness:}
\vspace{.3cm}
\begin{proposition}
\label{6.7}
We set $M=L\otimes_\Bbbk N$. Let $\Bbbk$ be an algebraically closed field, and suppose that
$\cd(J,N)=\ara(J,N)$. The following statements are equivalent:
 \begin{itemize}
\item[{(a)}]  $M$ is approximately Cohen--Macaulay with respect to $A \otimes_\Bbbk J$ which is not Cohen--Macaulay with respect to $A \otimes_\Bbbk J$;
\item[{(b)}] $M/(I \otimes_\Bbbk B)M$ is approximately Cohen--Macaulay with respect to $\mathcal{I}$ which is not Cohen-Macaulay with respect to $\mathcal{I}$.
 \end{itemize}
\end{proposition}
\begin{prf}
Let $M$ be approximately Cohen--Macaulay but not Cohen--Macaulay
with respect to $A \otimes_\Bbbk J$. By Corollary \ref{4.9}, this is
equivalent to saying that $N$ is approximately Cohen--Macaulay but
not Cohen--Macaulay with respect to $J$.
Notice that $M/(I \otimes_\Bbbk B)M\cong L/IL\otimes_\Bbbk N$ and $L/IL$ is Cohen--Macaulay with respect to $I$ with $\cd(I, L/IL)=\ara(I, L/IL)=0$.
 Thus, by Theorem \ref{4.7} this is equivalent to say that $M/(I \otimes_\Bbbk B)M$ is approximately Cohen--Macaulay but not Cohen--Macaulay with respect to $\mathcal{I}$.
\end{prf}

\vspace{.2cm}
\noindent{\bf Maximal depth property:}
\vspace{.3cm}

\begin{fact}
\label{6.8}{\em
Setting $I=0$ in Proposition \ref{5.1} and using
$\mcd(0, K)=0$ for every nonzero $R$-module $K$,
we obtain $\mcd(A\otimes_\Bbbk J, L\otimes_\Bbbk N)=\mcd(J,N)$. }
\end{fact}

\begin{proposition}
\label{6.9}
We set $M=L\otimes_\Bbbk N$. Then the following statements are equivalent:
\begin{itemize}
\item[(a)] $M$ has maximal depth with respect to $A\otimes_{\Bbbk}J$;
\item[(b)] $M/(I\otimes_{\Bbbk}B)M$ has maximal depth with respect to
$\mathcal{I}$.
\end{itemize}
\end{proposition}

\begin{prf}
Suppose first that $M$ has maximal depth with respect to
$A\otimes_{\Bbbk}J$. Then
$
\grade(A\otimes_{\Bbbk}J,M)
=
\mcd(A\otimes_{\Bbbk}J,M).
$
By Fact \ref{6.1} and Fact \ref{6.8},  it follows that
$
\grade(J,N)=\mcd(J,N).
$
Note that $M/(I\otimes_{\Bbbk}B)M
\cong
(L/IL)\otimes_{\Bbbk}N$. Moreover, $\grade(I,L/IL)=\mcd(I,L/IL)=0$ because $\cd(I,L/IL)=0$, see equation~(5.0.1). Thus, we obtain
\[
\begin{aligned}
\grade\left(
\mathcal{I},
M/(I\otimes_{\Bbbk}B)M
\right)
&=
\grade(I,L/IL)+\grade(J,N)\\
&=
\grade(J,N)\\
&=
\mcd(J,N)\\
&=
\mcd(I,L/IL)+\mcd(J,N)\\
&=
\mcd\left(
\mathcal{I},
M/(I\otimes_{\Bbbk}B)M
\right).
\end{aligned}
\]
The last equality follows from Proposition \ref{5.1}.
Thus, $M/(I\otimes_{\Bbbk}B)M$ has maximal depth with respect to
$\mathcal{I}$.

Conversely, suppose that
$M/(I\otimes_{\Bbbk}B)M$ has maximal depth with respect to $\mathcal{I}$.
Then
\[
\begin{aligned}
\grade\left(
\mathcal{I},
M/(I\otimes_{\Bbbk}B)M
\right)
&=
\mcd\left(
\mathcal{I},
M/(I\otimes_{\Bbbk}B)M
\right).
\end{aligned}
\]
Using
$
\grade(I,L/IL)=\mcd(I,L/IL)=0,
$
we obtain
$
\grade(J,N)=\mcd(J,N).
$
Hence
\[
\grade(A\otimes_{\Bbbk}J,M)
=
\mcd(A\otimes_{\Bbbk}J,M),
\]
and therefore $M$ has maximal depth with respect to
$A\otimes_{\Bbbk}J$.
Thus, (a) and (b) are equivalent.
\end{prf}

\vspace{.2cm}
\noindent{\bf Gorensteinness:}
\vspace{.3cm}

\begin{remark}{\em
\label{6.10}
Let $\cd(I, A)>0$ and $IL\neq L$.  We set $M=L\otimes_{\Bbbk}N$. In view of the discussion in this section, we may ask whether
$
M/(I\otimes_{\Bbbk}B)M
$
is Gorenstein with respect to $\mathcal{I}$ whenever $M$ is Gorenstein with respect to
$A\otimes_{\Bbbk}J$. We show that, in general, the answer is negative.

Indeed, suppose that $M$ is Gorenstein with respect to
$A\otimes_{\Bbbk}J$. Then
\[
\Ext^i_{A\otimes_{\Bbbk}B}
\left(
(A\otimes_{\Bbbk}B)/(A\otimes_{\Bbbk}J),
M
\right)\cong
L\otimes_{\Bbbk}\Ext^i_B(B/J,N)=0
\]
for all $i\neq\cd(A\otimes_{\Bbbk}J, A\otimes_{\Bbbk}B)=\cd(J,B)$.
Moreover, $L\otimes_{\Bbbk}\Ext^{\cd(J,B)}_B(B/J,N)\neq 0$. In particular, $\Ext^{\cd(J,B)}_B(B/J,N)\neq 0$.
On the other hand, $
M/(I\otimes_{\Bbbk}B)M
\cong (L/IL)\otimes_{\Bbbk}N
$.
By a consequence of the K\"unneth formula, we have
\[
\Ext^i_{A\otimes_{\Bbbk}B}
\left(
(A\otimes_{\Bbbk}B)/\mathcal{I},
M/(I\otimes_{\Bbbk}B)M
\right)
\cong
\bigoplus_{p+q=i}
\Ext^p_A(A/I,L/IL)
\otimes_{\Bbbk}
\Ext^q_B(B/J,N).
\]
Since $L/IL$ is an $A/I$-module, we have
\[
\Ext^p_A(A/I,L/IL)
\cong
\Ext^p_{A/I}(A/I,L/IL)=0
\quad\text{for all }p>0,
\]
while
$
\Ext^0_A(A/I,L/IL)
\cong L/IL.
$
Consequently,
\[
\Ext^i_{A\otimes_{\Bbbk}B}
\left(
(A\otimes_{\Bbbk}B)/\mathcal{I},
M/(I\otimes_{\Bbbk}B)M
\right)
\cong
(L/IL)\otimes_{\Bbbk}
\Ext^i_B(B/J,N).
\]
In particular,
\[
\Ext^{\cd(J, B)}_{A\otimes_{\Bbbk}B}
\left(
(A\otimes_{\Bbbk}B)/\mathcal{I},
M/(I\otimes_{\Bbbk}B)M
\right)
\cong
(L/IL)\otimes_{\Bbbk}
\Ext^{\cd(J, B)}_B(B/J, N)\neq 0.
\]
Therefore, $M/(I\otimes_{\Bbbk}B)M$ can not be Gorenstein with respect to $\mathcal{I}$. Indeed, if
$M/(I\otimes_{\Bbbk}B)M$ were Gorenstein with respect to $\mathcal{I}$,
then, by the definition of Gorensteinness with respect to an ideal, we would have
\[
\Ext^i_{A\otimes_{\Bbbk}B}
\left(
(A\otimes_{\Bbbk}B)/\mathcal{I},
M/(I\otimes_{\Bbbk}B)M
\right)=0
\quad\text{for all}\quad
i\neq
\cd(\mathcal{I},A\otimes_{\Bbbk}B)
=
\cd(I,A)+\cd(J,B),
\]
whereas the above computation shows that this Ext module is nonzero for $i=\cd(J, B)$ and $\cd(J,B)<\cd(I,A)+\cd(J,B)$, a contradiction.
Thus, the question in \cite{R5} or \cite{NR} has a negative answer for the Gorenstein property.}
\end{remark}

\vspace{.2cm}
\noindent{\bf Maximal Cohen--Macaulayness:}
\vspace{.3cm}

\begin{remark}
\label{6.11}{\em
Assume that $\cd(I,A)>0$ and $IL\neq L$. We set $M=L\otimes_{\Bbbk}N$.
In view of the preceding discussion, it is natural to ask whether
$
M/(I\otimes_{\Bbbk}B)M
$
is maximal Cohen--Macaulay with respect to $\mathcal{I}$ whenever
$M$ is maximal Cohen--Macaulay with respect to $A\otimes_{\Bbbk}J$.
In general, the answer is negative.

Indeed, suppose that $M$ is maximal
Cohen--Macaulay with respect to $A\otimes_{\Bbbk}J$. Then
\[
\grade(A\otimes_{\Bbbk}J,M)
=
\cd(A\otimes_{\Bbbk}J,A\otimes_{\Bbbk}B).
\]
By Fact \ref{6.1}, it follows that $\grade(J,N)=\cd(J,B)$.
On the other hand,
$
M/(I\otimes_{\Bbbk}B)M
\cong (L/IL)\otimes_{\Bbbk}N.
$
Hence,
\[
\begin{aligned}
\grade\left(
\mathcal{I},
M/(I\otimes_{\Bbbk}B)M
\right)&=\grade(I,L/IL)+\grade(J,N)\\
        &= \grade(I,L/IL)+\cd(J,B).
\end{aligned}
\]
Since $IL\neq L$, we have $L/IL\neq 0$. Moreover, $\cd(I,L/IL)=0$, and therefore $\grade(I,L/IL)=0$.
Consequently,
\[
\grade\left(
\mathcal{I},
M/(I\otimes_{\Bbbk}B)M
\right)
=\cd(J,B).
\]
However,  $\cd(\mathcal{I},A\otimes_{\Bbbk}B)
=\cd(I,A)+\cd(J,B).
$
Since $\cd(I,A)>0$, we obtain
$\grade\left(
\mathcal{I},
M/(I\otimes_{\Bbbk}B)M
\right)\neq
\cd(\mathcal{I},A\otimes_{\Bbbk}B)$.
Thus, $M/(I\otimes_{\Bbbk}B)M$ is not maximal
Cohen--Macaulay with respect to $\mathcal{I}$.
Consequently, the corresponding assertion in
 \cite{R5} or \cite{NR}
does not hold in general.}
\end{remark}

\end{document}